\documentclass[twoside,10.5pt]{article}
\usepackage{mathrsfs}
\usepackage{pifont}
\usepackage{amsmath}
\usepackage{amsthm}
\usepackage{txfonts}
\usepackage{geometry}
\usepackage{latexsym}
\usepackage{amssymb}
\usepackage{graphicx}
\usepackage{geometry}
\usepackage{xcolor} 
\begin{document}
\begin{center}
{\LARGE{Probabilistic Analysis of Balancing Scores for Causal Inference}}\\[20pt]
Priyantha Wijayatunga$^1$
\end{center}
$^1$ Department of Statistics, Ume\r{a} School of Business and Economics, Ume\r{a} University, Ume\r{a}, Sweden.

Correspondence: Priyantha Wijayatunga, Department of Statistics, Ume\r{a} School of Business and Economics, Ume\r{a} University, Ume\r{a} SE-901 87,  Sweden. Tel: ++46-90-786-6610. E-mail: priyantha.wijayatunga@stat.umu.se
\\[15pt]
Received: December 28, 2014 \quad  Accepted: January 27, 2015 \quad Online Published: November XX, 201X \par
doi:10.5539/jmr.v7n2px  \quad \quad \quad    URL: http://dx.doi.org/10.5539/jmr.v7n2px \\[10pt]


\textbf{Abstract}

Propensity scores are often used for stratification of treatment and control groups of subjects in observational data to remove confounding bias when estimating of  causal effect of the treatment on an outcome in so-called potential outcome causal modeling framework. In this article, we try to get some insights into basic behavior of  the propensity scores in a probabilistic sense. We do a simple analysis of their usage confining to the case of discrete confounding covariates and outcomes. While making clear about behavior of the propensity score our analysis shows how the so-called prognostic score can be derived simultaneously. However the prognostic score is derived in a limited sense in the current literature whereas our derivation is more general and shows all possibilities of having the score. And we call it outcome score. We argue that application of both the propensity score and the outcome score is the most efficient way for  reduction of dimension in the confounding covariates as opposed to current belief that the propensity score alone is the most efficient way.    


\textbf{1. Introduction}

In many real world situations whether it is socio-economic, medical, etc. estimation of causal effect of some event or treatment on a certain outcome is of interest. Taking a certain medicine for a certain disease by patients and participating in a certain job training program by unemployed individuals in order to find employment are two examples. Sometimes it may be unethical or infeasible to assign randomly each subject either to the treatment or the control group in order to do a randomized study. However, for policy makers, medical professionals, etc.,  evaluation of the causal effects of those treatments may be of interest. In the absence of the random assignment of the treatment which  is considered to be the gold standard for calculation of the causal effect of it they may only have observed data on a collection of subjects  who either have taken the treatment or not.  When  the causal effect of the treatment on the outcome needs to be identified from such observational data it needs to control for a 'sufficient' subset of all the confounders of the causal relation (see, for example Vansteelandt et al. 2010). That is, direct comparison of the whole treatment group and the control group as in the case of randomized treatment assignments should not be done.  Subgroups of subjects with same confounder values (of the sufficient set) in the treatment group and in the control group should be compared first. This is called stratification of the data sample. And then those resulting comparisons should be weighted averaged over the subgroups  to yield the total comparison where the weights are the proportions of sizes of the respective subgroups.  Note that the confounders are factors affecting the subjects to take the treatment or not and simultaneously affecting the outcomes of the subjects in some way. Therefore the effect of the treatment on the outcome can be confounded with the effects of these confounding factors when they are present. However finding all the confounders of a given causal relation is a difficult task.   

And sometimes controlling for the confounders can be difficult, for example, when the confounders are high dimensional. Then it may be difficult to find the treatment and the control subgroups of the subjects of sufficient sizes with same confounder values.   A popular way to increase the sizes of the treatment and the control subgroups that should be compared, ideally when all the confounders are identified, is to  use so-called propensity scores (Rosenbaum \& Rubin 1983). The propensity score is the  conditional probability of receiving the treatment given the values of all those  observed pre-treatment confounding covariates. They are used in the causal inference method  of potential outcome modeling framework (Rubin 1974, 2007, 2008 and Holland 1986) for finding comparable subgroups of the treated subjects and untreated (controlled) subjects, i.e., to have a new stratification on the data sample from the old stratification. Note that when data sample is stratified with respect to assumed confounders then the number of strata is that of all possible configurations of the confounders. But if there exist some confounder configurations where their respective probabilities of being treated are equal to each other (when the propensity scores related to those configurations/strata are equal) then a new stratification can be done by merging those strata together. If it is done then it results in lesser number of strata in the data sample since the number of distinct propensity score values is lesser than the number of (distinct) configurations in the assumed confounders.  Therefore, generally the numbers of subjects in each strata of new stratification are larger  than or equal to those of the old stratification for the data sample. However when a set of potential confounders are available the potential outcome causal modeling framework offers no clear way to select the ones that should be included in the propensity score calculation.  But this selection affects the new stratification and so may be the causal effect estimates. 

Since propensity scores are conditional probabilities, they  are usually estimated by applying a logistic regression model that uses all the potential confounders, especially when the data are sparse.  However generally it is not known the true propensity score model in a given situation. Though in some places in literature (for example, Shah et al. 2005 and St\"umer et al. 2006)  it has been shown that the propensity score based causal effect estimation methods and traditional methods that are based on linear, logistic and Cox regression generally yield similar results, in some other places (for example, Martens et al. 2008), it has been demonstrated  that the propensity score based methods yield causal effect estimates that are generally closer to the true values than estimates based on the regression based methods. Note that in regression based methods adjustment for the confounders is done by applying the regression models on the outcome with the treatment variable and the confounders. 

Here, we try to get some insights into the basic behavior of  the propensity scores and balancing scores in a simple probabilistic sense where the balancing scores are more general scores that can be used to increase the sizes of the treatment and the control subgroups that should be compared for the estimation of the causal effect of the treatment. The propensity score is a special case of the balancing score. Though the propensity scores are used in many ways, here we confine to a fundamental analysis of their simple usage. Our focus is on their behavior in the context of causal inference in a probabilistic sense.  Therefore we assume that we have a 'sufficient' subset of all potential confounders, that are able to remove the confounding bias. We are not addressing the problems of unobserved confounders. Furthermore we do not address any selection bias problem often encountered in the estimation of causal effects.  

This simple analysis shows us how to derive so-called prognostic score (Hansen 2008) with similar arguments. However the prognostic score found in current literature is derived in a limited sense. Our derivation of this score that we call potential outcome score is more general and furthermore we look at all possibilities in doing so in a probabilistic sense. Therefore, it can be clearer to non-statisticians such as computer scientists, social scientists, etc. whereas Hansen's is more through a statistical analysis. In current literature (Hansen 2008) it has been described other forms of prognostic scores and the reader is referred to them for more details on them.  Importantly, our probabilistic analysis derives both the propensity score and the outcome score simultaneously. Usage of the outcome (prognostic) score is fundamentally same as that of the propensity score and we show that having two similar balancing scores allows us to do stratification of the data sample in both ways but one at a time stepwise. Furthermore we show that by using both scores effectively one can possibly obtain the smallest number of strata for the data sample.  Currently the propensity score is extensively used in practice but the prognostic score is not used so often, perhaps due the limited sense derivation of it. So, we think that our derivation will encourage the practitioner using the outcome score more often. We confine to  simple case of discrete confounding covariates and binary  treatment and the outcome variables and assume sufficiently large sample of data for probability estimations and calculations. 
\vspace{0.3cm}

\textbf{2. Observational Studies}

We consider the simple situation where one is interested in evaluating the effect of some exposure or treatment on a certain outcome that can either be a success or a failure. Let us denote the  treatment by a binary variable $Z$ where $Z=1$ when the treatment is implemented and $Z=0$ when it is not, and the outcome by a binary variable $Y$ where $Y=1$ when a success is observed and $Y=0$ when a failure is observed for each subject concerned. In the potential outcome framework for causal inference (Holland 1986) it is accepted existence of a pair of  potential outcome variables, say, $(Y_1, Y_0)$ where $Y_i$ is the outcome that would have been observed  had the treatment $Z=i$ for $i=1,0$. Note that then the observable outcome $Y$ satisfies  $Y=ZY_1 + (1-Z)Y_0$. Then a randomized experiment is when the potential outcomes are independent of the treatment assignment, written as $(Y_0,Y_1) \perp Z $; each subject receives the treatment without considering its future outcome. But in observational data this may not hold because subjects do not receive the treatment independent of their future outcomes, therefore characteristics of the subjects in the treatment group may differ from those of the control group systematically. This is a situation where the treatment effect is confounded with some external factors, i.e.,  the treatment and the outcome are confounded. Therefore the treatment and the control groups cannot be compared directly to evaluate the effect of the treatment as in case of randomized study. Then the assumption is modified such that the potential outcomes are conditionally independent of the treatment assignment given all those pre-treatment confounding factors. When all such confounders are denoted by multivariate variable $X$  then the assumption is written as $(Y_0,Y_1) \perp Z \vert X$ and it is called 'assumption of no unmeasured confounders'. Here we assume that $X$ is discrete, therefore any continuous covariate is discretized.  That is, at each stratum of $X$ the treatment assignments are assumed to be randomized.  For this mimicking  of the randomization of the treatment assignment within each stratum of $X$, firstly one should have found all the confounders $X$. However this assumption cannot be tested even if all the potential confounders are found.  

From a given set of directly or indirectly related pre-treatment confounders of the treatment and the outcome finding a set of confounders that should be controlled for is somewhat problematic and the potential outcome framework offers no clear way to do it. Note that it may be redundant to confrol for all such confounders or sometimes we need to omit some of them since inclusion of them can cause extra bias (Shrier 2008, Rubin 2009, Sj\"{o}lander 2009 and Pearl 2009a). However the causal graphical modeling framework of Pearl and his colleagues (Pearl 2009 and references therein) offers one called 'back door criterion' to choose a set of confounders in order to identify the causal effect of $Z$ on $Y$, i.e., to estimate without bias.  When a graphical model is done with $Z$ and $Y$ and ideally all their causal factors both direct and indirect the criterion can find a sufficient set of confounders on which one should control for the estimation of the causal effect. The set is called 'admissible' or 'deconfounding' set. And considering some covariates as confounders by ignoring such criterion can introduce  further bias (p. 351 of Pearl 2009).   In our analysis we confine to the case of that $X$ is an admissible set of confounders and stated otherwise. 

Furthermore, sometimes observed data may suffer from selection bias problem, i.e., data are selected according to some selection criteria and they do not represent the target population of study (Spirtes et al 2001  and references therein). It is discussed in literature (Cai \& Kuroki 2008) how to work out Pearl's graphical modeling criteria for selecting a sufficient set of confounders in the presence of selection bias. In our analysis we assume that our data sample is not suffering from the selection bias problem.

In what follows we use Rubin's causal model (Holland 1986). Let us define the individual causal effect for an individual, say, $j$ with $X=x$ as $\tau^{j}(x)=Y_1^j -Y_0^j$. But it is clear that no subject has both the values of $Y_1$ and $Y_0$ observed therefore we cannot have $\tau^j(x)$ numerically. So we need a mechanism to get it but it is right at our hands; the randomization of treatment within  each stratum of $X$ which is an admissible set, the assumption of no unmeasured confounders. This is also called assumption of strong ignorable treatment assignment  (Rosenbaum \& Rubin 1983). That is, within each stratum $X=x$ if we know a subject is treated ($Z=1$) therefore we observe $Y_1=Y$ but $Y_0$ is not known. The latter can be known by any other subject in the stratum who is not treated ($Z=0$); two quantities are conditionally exchangeable.  Here the word conditionally is to mean within the stratum. Note that in a randomized study these two quantities are marginally exchangeable. And similarly for any subject that is not treated ($Z=0$). Therefore, as if the observed data are from randomization within each level of $X$, we can calculate the  average causal effect for the subpopulation of all individuals with $X=x$, say, $\tau(x)$ by 
\begin{align}
\tau(x) & =  E[Y_1 \vert X=x] -  E[Y_0 \vert X=x] =   E[Y_1 \vert X=x,Z=1] -  E[Y_0 \vert X=x,Z=0] \nonumber \\
 \tau(x)  & =   E[Y \vert X=x,Z=1] -  E[Y \vert X=x,Z=0]  \label{taux} 
\end{align}
where the expectation $E$ should be taken over whole subpopulation with $X=x$. Since this mechanism applies for all the strata of $X$, we can calculate  the average causal effect for the whole population, say, $\tau$    as follows.
\begin{align}
\tau &= E [ E[Y \vert Z=1, X=x] -  E[Y \vert Z=0, X=x]] \nonumber \\
\tau  &= \sum_y y \sum_i p(Y=y \vert X=x_i,Z=1)p(X=x_i) - \sum_y y \sum_i p(Y=y \vert X=x_i,Z=0)p(X=x_i)  \label{tau}
\end{align}

The above estimate for $\tau$ is analytically equal to that we get by the estimation of the causal effects using interventions in the causal graphical models, also called do-calculus, (Pearl 2009 and Lauritzen \& Richardson 2002),  therefore two frameworks are equivalent in this case.  In literature it is discussed  (Wijayatunga 2014) about this relation considering how some of the causal effect estimators found in the potential outcome causal modeling framework are related to causal graphical model application.
\small
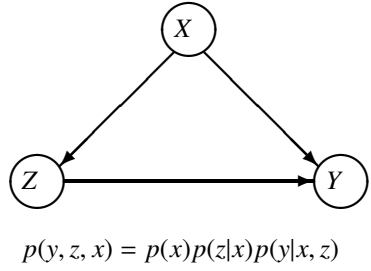
\begin{figure} 
\begin{center}
\setlength{\unitlength}{1mm}
\begin{picture}(100,40)(-50,-10)
\thicklines
\put(0,20){\circle{7}}
\put(-2,19){$X$}

\put(-20,0){\circle{7}}
\put(-22,-1){$Z$}

\put(20,0){\circle{7}}
\put(18,-1){$Y$}

 \put(-2,17){\vector(-1,-1){15}}
 \put(2,17){\vector(1,-1){15}}
 \put(-16.5,0){\vector(1,0){33}}

\put(-22, -10){$p(y,z,x)=p(x)p(z\vert x)p(y\vert x,z)$}
\end{picture}
\end{center}
\caption{ \label{simple.bn}  Causal graphical model for causal relations among $X,Z$ and $Y$.}
\end{figure}
\normalsize 
The causal relationships between $X$, $Y$ and $Z$ in our context can be represented as a causal network model $p(y,z,x)=p(x)p(z\vert x)p(y\vert x,z)$ as shown in the Figure \ref{simple.bn} where an arrow represents a causal relationship between a cause and an effect such that the arrow starts from the cause and pointed to the effect. If we intervene on $Z$ to have its value $z$, written as $do(Z=z)$ for $z=0,1$, then we have the intervention distribution of $Y$ and $X$, 
\begin{equation}
p(Y= y, X=x \vert do(Z=z)=\frac{p(X=x)p(Z=z \vert X=x)p(Y=y \vert Z=z,X=x)}{p(Z=z \vert X=x)}  \label{pyx_z}
\end{equation}
So we have $p(Y=y \vert do(Z=z)) = \sum_i p(Y=y \vert Z=z,X=x_i)p(X=x_i)$. The causal effect of the treatment option $Z=1$ compared to the control option $Z=0$, say $\rho$, is defined as the difference between the  expectations of $Y$ under two interventions,
\begin{align}
  \rho &= \sum_y  y p(Y=y \vert do(Z=1))  - \sum_y y p(Y=y \vert  do(Z=0))  \nonumber \\
  \rho     &=  \sum_y y \sum_i p(Y=y \vert Z=1,X=x_i)p(X=x_i) -  \sum_y y \sum_i p(Y=y \vert Z=0,X=x_i)p(X=x_i)    \label{rho}
\end{align}
which shows that $\tau$ is equal to $\rho$, so we have seen that the  strong ignorable treatment assignment assumption  in the potential outcome model is equivalent to implementing intervention operations in the probability distributions when in both cases the confounding factors are the same, i.e., they yield analytically the same causal effect estimates. In fact one can see that the potential outcome distribution under hypothetical treatment assignment assuming no further unmeasured confounders  and the outcome distribution under the same intervention on the treatment with the same set of confounders are the same. For $i=0,1$, we have that
\begin{align}
p(Y_i=y) &= \sum_k p(Y_i=y \vert x_k)p(x_k) = \sum_k p(Y_i=y \vert Z,x_k)p(x_k) \nonumber \\
 &= \sum_k p(Y=y \vert Z=i,x_k)p(x_k)  =p(Y=y \vert do(Z=i)) \label{pYi}
\end{align}
More generally if the treatment $Z$ affects outcome $Y$ directly or indirectly through a mediator $Z'$ or in both ways, then we can show that Eq. (\ref{pYi}) holds, therefore in those cases two causal effect estimates are analytically the same.

Since by definition all the possible values of $X$ are fixed, there should be enough data in the contingency table of $X,Z$ and $Y$ for accurate estimation of above casual effect. Sometimes some of the strata can be merged together, often with respective to the variable $X$, for example, to avoid data sparsity. This results in a new stratification for the confounders thus forming a new variable, say, $X'$.  Since it is accepted that randomization is the gold standard for causal inference,  we try to create a similar situation for performing causal inference with the observational data. In a randomized study the probability of the assignment of the treatment is the same for all subjects who are characterized by pre-treatment information/covariate $X$, i.e., $p(Z=1 \vert X=x)$ is the same for all $x$. Therefore all the strata of $X$ are  pooled together both  in the treatment group and the control group  separately for subsequent comparison the two for identifying the causal effect of the treatment. Mimicking this within the observational data means whenever we have $p(Z=1 \vert X=x_i)=p(Z=1 \vert X=x_j)$, (say, it is $p$) for some $i \neq j$ we can pool the two strata $X=x_i$ and $X=x_j$  together for strata-wise comparison of treatment and non-treatment groups. This is a consequence of simple fact that when above case is true then $p(Z=1 \vert X=x_i \textrm{ or } X=x_j)=p$. So, by this way, we can reduce the dimension of the conditioning set $X$.  Note that if we do so, then we have a new conditioning (confounding) variable $X'$ with lesser number of configurations than that of the original confounding variable.

 To see that we get same causal effect estimates in both the cases of before strata merging and after strata merging, without loss of generality, take $i=1$ and $j=2$ for the observed data sample and then consider the first two terms in of the summation, $ \sum_k p(Y=y \vert X=x_k,Z=1)p(X=x_k)$  is
\begin{align*}
 \sum_{k=1}^2 & p(Y=y \vert X=x_k,Z=1)p(X=x_k)  \\
   &= \sum_{k=1}^2  p(Y=y \vert X=x_k,Z=1)p(X=x_k)  \frac{p(Z=1 \vert X=x_k)}{p(Z=1 \vert (X=x_
1 \textrm{ or } X=x_2) )} \\
  &= \sum_{k=1}^2 p(Y=y \vert X=x_k,Z=1) \frac{p(Z=1, X=x_k)}{p(Z=1, (X=x_1 \textrm{ or } X=x_2 ))} p(X=x_1 \textrm{ or } X=x_2 ) \\
&=  \frac{p(X=x_1 \textrm{ or } X=x_2 )}{p(Z=1, (X=x_1 \textrm{ or } X=x_2) )}  \sum_{k=1}^2 p(Y=y, X=x_k,Z=1) \\
&= p(Y=y \vert Z=1, (X=x_1 \textrm{ or } X=x_2 ))  p(X=x_1 \textrm{ or } X=x_2 )
\end{align*}
Therefore, we define a new variable, $X'$ from $X$ as follows; $X'=x_1'$ when $X=x_1$ or $X=x_2$ and $X'=x_{i-1}'$ when $X=x_i$ for $i=3,...,n$. Then, from Eq. (\ref{tau}) we have
\begin{align*}
\tau &= \sum_y y \Bigg\lbrace   p(Y=y \vert Z=1, (X=x_1 \textrm{ or } X=x_2))  p(X=x_1 \textrm{ or } X=x_2) +   \sum_{k=3}^n p(Y=y \vert X=x_k, Z=1)p(X=x_k)  \Bigg\rbrace \\
 & \qquad{} - \sum_y y \Bigg\lbrace   p(Y=y \vert Z=0, (X=x_1 \textrm{ or } X=x_2))  p(X=x_1 \textrm{ or } X=x_2)  +   \sum_{k=3}^n p(Y=y \vert X=x_k, Z=0)p(X=x_k)  \Bigg\rbrace \\
 &= \sum_y y \sum_{k=1}^{n-1} p(Y=y \vert X'=x'_k,Z=1)p(X'=x'_k) - \sum_y y \sum_{k=1}^{n-1} P(Y=y \vert X'=x'_k,Z=0)p(X'=x'_k) 
\end{align*}
That is, in general, as long as there are equalities in the conditional probabilities $p(Z=1 \vert X=x)$ with respect to different values of $X$, those strata can be merged together into one, thus forming a new conditioning variable $X'$ out of $X$. Then the causal effect can be estimated with respect to this new confounding variable and arrive at the same value as in the case of original confounding variable $X$. And with this new variable one obtains some estimation power (more data for some of the value categories of confounding variable). In the case of existence of such equalities there exists a mapping which is many-to-one and onto from the set of values (state space) of $X$ to that of $X'$; for each value $x'$ of $X'$ there exists at least one value $x$ of $X$ such that $p(Z=1 \vert X'=x')=p(Z=1 \vert X=x)$. Therefore $X'$ takes lesser number of values than that of $X$. However in the potential outcome causal model the real use of $X'$ is implicit. Instead for each confounder value $x$ the propensity score, the probability that a subject with its  confounder value $x$ taking the treatment $p(Z=1 \vert X=x)$ is used explicitly in causal effect calculations. Though the statisticians may be aware, due to this implicitness applied users of statistics may interpret that the propensity scores generally solve their causal inference problems. It should be understood that by such mergers of strata of the confounder set one obtain same causal effect estimates as in the case of without such mergers under strict equalities among propensity scores related to the mergers. However one effective use of propensity score is when one considers approximate equalities among them for such mergers. And in those contexts it may be possible that two causal effect estimates to differ numerically to some extent where difference may depends on the closeness of the approximations considered. 
\vspace{0.3cm}

\textbf{3. Balancing Scores}

Here we will have simple look at conditional independence relations related to the balancing scores and in particular the propensity scores. And this allows us to see that there exists another score that is similar to the propensity score, which we call potential outcome score. However a same kind of score that is called prognostic score was derived in literature (Hansen 2008) but here our concern is on simultaneous derivation of all possible scores probabilistically. The effects of the outcome score in the causal effect calculation seems to be generally similar to that of the propensity score. However, even though such a score is obvious, most of the discussions in the literature so far are around the propensity scores.   Note that the propensity scores can reduce the dimension of the conditioning variable set maximally, i.e., can obtain $X'$ from $X$ such that $X'$ has smallest number of values, when using $Z$ and $X$ only. But with the outcome score explained below, which is a function of $Z, X$ and $Y$ we can also reduce the dimension of the conditioning set irrespective of any other score defined with $X$ and $Z$ such as the propensity score. That is,  the outcome score does not prohibit using any other score, for example, the propensity score in an earlier stage or later stage of the dimension reduction process of the conditioning set $X$. And so does the propensity score. Therefore we show that maximal ('more efficient' in sense of Hansen 2006) dimension reduction is  obtained through the application of both the propensity score and the outcome score, ideally one at a time. 

First we look at how  the propensity score works (which we already have seen) but little more formerly and then analogously we define the outcome score. Let $X$ be a discrete covariate vector with a positive probability, i.e., for all possible values $x$ of $X$, $p(X=x)>0$ and $Z$ is a binary random variable taking values in $\lbrace 0, 1 \rbrace$ with positive probability.  Let $b$ be an onto function of $X$ (co-domain is the image of $b$ under that mapping), therefore  $b(X)$ is a random variable with a positive probability and the cardinality of the state space of $b(X)$ is at most that of $X$. Recall that the cardinality of a set is the number of distinct values in it. As earlier, let us denote that the conditional independence of $X$ and $Y$ given $Z$, that is, when $P(X,Y \vert Z)=P(X \vert Z)P(Y \vert Z)$ by $X \perp Y \vert Z$ and marginal independence of $X$ and $Y$ that is, when $P(X,Y )=P(X)P(Y )$ by $X \perp Y$.

\emph{3.1 Propensity Scores}

A balancing score (Rosenbaum \& Rubin 1983), say, $b(X)$, is a function of $X$ such that $X \perp Z \vert b(X)$ and then the trivial balancing score is $X$ itself. Note that $X \perp Z \vert X$ is true for any $X$ and $Z$ where $X$ has a positive probability. Furthermore it is stated that  more useful balancing scores are  many-to-one functions of $X$.  To see this, without loss of generality,  suppose a many-to-one and onto function $b$ such that $b(x_1)=b_1=b_2=b(x_2)$ and $b(x_i)=b_i$ for $i=3,...,n$ where $b_3,...,b_n$ are $n-2$ number of distinct values that are not equal to $b_1$. That is, when state space of $X$ is $\lbrace x_1,x_2,....x_n \rbrace$  then that of $b(X)$ is $\lbrace b_1=b_2,b_3...,b_n\rbrace $, so $b$ is a many-to-one function of $X$. Let $0< p(X=x_i,Z=z)<1$ for $i=1,...,n$ and $Z=0,1$, therefore  so do their marginals. Then
\begin{align*}
p(X = x_1, Z=z \vert b(X)=b_1) & = \frac{p(X=x_1,Z=z,b(X)=b_1)}{p(b(X)=b_1)} \\
  &= \frac{p(X=x_1,Z=z,(X=x_1 \textrm{ or } X=x_2))}{p(b(X)=b_1 )}  \\
  &= \frac{p(X=x_1,Z=z)}{p(b(X)=b_1 )} \frac{p(Z=z,b(X)=b_1)} {p(Z=z,b(X)=b_1)} \\
  & =    p(Z=z \vert b(X)=b_1)  p(X=x_1 \vert b(X)=b_1)   \times     \frac{p(Z=z \vert X=x_1)}{p(Z=z \vert b(X)=b_1)} 
\end{align*}
Therefore  $X \perp Z \vert b(X) =b_1$ if and only if $p(Z=z \vert X=x_1)=p(Z=z \vert X=x_2)$. And since $b_1=b_2$ this is true for $b(X)=b_2$.  Similarly we get the required results when  $X=x_2$. Since $x_i$ is mapped to $b_i$ for each $i=3,...,n$ the relation $X \perp Z \vert b(X)=b_i$   hold for those $b_i$'s too. Therefore $X \perp Z \vert b(X)$. Note that this result is true for any many-to-one mapping $b$ of $X$ such that $b_i=b(x_i)=b(x_j)=b_j$ whenever $p(Z=z \vert X=x_i)=p(Z=z \vert X=x_j)$ for $i \neq j$. Therefore we have the following theorem.

\textbf{Theorem 1}~\emph{When we have that $X \perp Z \vert b(X)$ where $b$ is a function of $X$ then the cardinality of the state space of  $b(X)$ is lesser than that of  $X$ by $m-1$ whenever there are $m$ distinct values of $X$ such that $p(z \vert x_{i_1} )=...=p(z \vert x_{i_m})$ for $1\leq i_1 < ... < i_m \leq n$.}

Proof is clear from above. Note that this theorem applies sequentially as long as there are many distinct such equality sequences. At each step one can have a new conditioning variable $X'\equiv b(X) $ from the old one $X$. No further dimension reduction is possible unless considering approximate equalities. If we let $b(X)$ to be the propensity score $p(Z=1 \vert X)$ then $b(X)$ has the maximal dimension reduction, i.e., the propensity score defines the conditioning variable with minimal cardinality (of its state space) for conditional independence of $Z$ and $X$.  Note that no other function that has lesser cardinality in its image can make $Z$ and $X$ conditional independent given it, since any such function should be built considering $ p(Z=1 \vert X)$. 

It is clear from above that if we have another $b^*(X)=p^*(Z=1 \vert X)$ such that $b^*(X) \not\equiv b(X)$ but there is a one-to-one relationship between $b(X)$ and $b^*(X)$ such that whenever $b(x_i)=b(x_j)$ then $b^*(x_i)=b^*(x_j)$ for $i \neq j$ then we have $X \perp Z \vert b^*(X)$. Let $b(X)$ is the true propensity score that gives us $X'$ and $b^*(X)$ is a misspecified propensity score that gives us $X^*$ then, still we obtain the correct subgroups to compare in order to evaluate causal effect since $X'$ and $X^*$ should coincide. And we can obtain correct causal effect estimates as long as  we have $P(Y \vert X',Z)=P(Y \vert X^*,Z)$ and $P(X')=P(X^*)$. Note that we will have them naturally because only difference between the two cases is the naming of strata of confounders.  However  if we use $p(Z=1 \vert X)$ and $p^*(Z=1 \vert X)$ explicitly in the causal effect formula we will generally have different causal effect estimates for the two cases  since then Eg. (\ref{tau}) becomes 
\begin{equation*}
\tau =\sum_{y,x} y \frac{p(Y=y,Z=1,X=x)}{p(Z=1 \vert X=x)}-\sum_{y,x} y \frac{p(Y=y,Z=0,X=x)}{p(Z=0 \vert X=x)}
\end{equation*}
So, in the potential outcome causal model where propenisty scores are explicitly used misspecification of propensity scores may yield incorrect causal effect estimates. Furthermore the misspecification is such that, for some distinct $i,j$ we have $b^*(x_i) \approx b^*(x_j)$ when $b(x_i) \neq b(x_j)$ then we may end up pooling two strata $X=x_i$ and $X=x_j$ when we are using misspecified propensity score.

\emph{3.2 Searching for Other Scores}

Now for some $Y$ let $Y \perp Z \vert X$. Then $Y \perp Z \vert b(X)$ if and only if, without loss of generality, the function $b$ is defined as in the above case but with a revised condition that $p(Z \vert X=x_1)=p(Z \vert X=x_2)$ or $p(Y \vert X=x_1)=p(Y\vert X=x_2)$. 
\begin{align*}
p(Y,Z \vert  b(X)=b_1) & = \frac{p(Y,Z, b(X)=b_1)}{p( b(X)=b_1)} \\
     &= \frac{p(Y,Z, (X=x_1 \textrm{ or } X=x_2))}{p(X=x_1 \textrm{ or } X=x_2)} \\
      &= \frac{p(Y,Z, X=x_1)}{p(X=x_1 \textrm{ or } X=x_2)}  + \frac{p(Y,Z,  X=x_2)}{p(X=x_1 \textrm{ or } X=x_2)} \\
     & =  w p(Y,Z \vert  X=x_1) + (1- w) p(Y,Z \vert X=x_2) \\
     & =  w p(Y \vert  X=x_1)p(Z \vert X=x_1)   + (1- w) p(Y \vert X=x_2)p(Z \vert X=x_2) \\
      & =  p(Y \vert b(X)=b_1)p(Z \vert b(X)=b_1)
\end{align*}
where $ w=\frac{p(X=x_1 )}{ p((X=x_1 \textrm{ or } X=x_2))}$. Again note that this result is true for any many-to-one mapping $b$ of $X$ such that $b_i=b(x_i)=b(x_j)=b_j$ whenever $p(Z=z \vert X=x_i)=p(Z=z \vert X=x_j)$ or $p(Y=y \vert X=x_i)=p(Y=y \vert X=x_j)$ for $i \neq j$. Therefore we have the following theorem.

\textbf{Theorem 2}~\emph{Let $Y \perp Z \vert X$ then we have $Y \perp Z \vert b(X)$ where $b$ is a many-to-one function of $X$ such that the  cardinality of the state space of $b(X)$ is lesser than that of $X$ by $m-1$ whenever  there are $m$ distinct values of $X$ such that either $p(z \vert x_{i_1} )=...=p(z \vert x_{i_m})$ for $1\leq i_1 < ... < i_m \leq n$ or $p(y \vert x_{j_1} )=...=p(y \vert x_{j_m})$ for $1\leq j_1 < ... < j_m \leq n$.}

Again proof is clear from above. Similar to earlier theorem this theorem also applies sequentially as long as there are many distinct such equality sequences. At each step one can have a new conditioning variable $X'\equiv b(X) $ from the old one $X$. No further dimension reduction is possible unless considering approximate equalities. There can be some dimension reduction  when conditioning with a balancing score rather than with the original conditioning variable where the balancing score is defined appropriately. There are two ways to do it. That is,  when some equalities among probabilities $P(Z \vert X)$ or $P(Y \vert X)$ exist we can form smaller size contingency table on $b(X), Y$ and $Z$ with the data sample, that have more entries in its cells than the original contingency table of $X,Y$ and $Z$. Forming this new contingency table can be done in term of $Z$ (with $P(Z \vert X)$) and then in terms of $Y$ (with $P(Y \vert X)$) or vice versa. 

It should be noted   that when some operations with probabilities $P(Z \vert X)$ are done some additional possibilities may arise wih respect to probabilities $P(Y \vert X)$. To see this suppose, for example,  let we observe that $p(z \vert X=x_1)=p(z \vert X=x_2)$ and therefore we merge two strata $X=x_1$ and $X=x_2$, thus resulting a new conditioning variable, say,  $X'$($ \equiv b(X)$).   As previously, let $X'=x'_1$ when $X=x_1$ or $X=x_2$ and $X'=x'_{i-1}$ when $X=x_i$ for $i=3,...,n$ where set of distinct values of $X$ is $\{x_1,...,x_n\}$. But this operation results in that $p(y \vert X'=x'_1)=w p(y \vert X=x_1)+ (1-w) p(y \vert X=x_2)$. Note that $p(y \vert X=x_1) \leq p(y \vert X'=x'_1) \leq p(y \vert X=x_2)$ when $p(y \vert X=x_1) \leq p(y \vert X=x_2)$. Then additionally relations such as $p(y \vert X'=x'_1)=p(y \vert X'=x'_i)$ for some $i >1$ may arise. If so, one can do the dimension reduction with respective to $P(Y \vert X')$ in the contingency table of $X',Y$ and $Z$. 

The other kind of balancing score found in causal inference literature is so-called prognostic score (Hansen 2008). The prognostic score is derived using $P(Y_0 \vert X)$, i.e., it is defined as the predicted outcome under control condition. It is estimated by fitting a model of the outcome in the control group and then using that model to obtain predictions of the outcome under the control condition for all subjects (Stuart et al. 2013). However the above theorem can be applied to the potential outcome causal model where it is assumed $(Y_0,Y_1) \perp Z \vert X$. Then one can find $b$ satisfying $(Y_0,Y_1) \perp Z \vert b(X)$ where $b$ is defined such that $b(x_i)=b_i= b_j=b(x_j)$ whenever, for example, either (1) $p(z \vert x_i )=p(z \vert x_j)$ or (2) $p(Y_0=y_0 \vert x_i )=p(Y_0 =y_0 \vert x_j)$ and $p(Y_1=y_1 \vert x_i )=p(Y_1=y_1 \vert x_j)$, for $i \neq j$, $i,j=1, \dots, n$. The above case (1) refers to the propensity score and case (2) should refer to the prognostic score. However the score derived in case (2) is a more general one than the prognostic score that is found in current literature. Therefore, we call this score the potential outcome score that is discussed in more details later.   Though the propensity and the prognostic scores are derived quite distinctively in the past, our simple probabilistic analysis of the conditional independence relation $Y_0,Y_1 \perp Z \vert X$ derives the two scores at once. And furthermore, by this way it may be easy to see how the two  scores are related to each other and also to use them simultaneously.

\emph{3.3 Potential Outcome Scores}

Now we can look at how to obtain balancing scores from a construction of $b$ through equality relations with  conditional probabilities $P(Y \vert X)$. Consider the case where $P(Y_0 \vert X=x_i)=P(Y_0 \vert X=x_j)$ and $P(Y_1  \vert X=x_i)=P(Y_1 \vert X=x_j)$ for some $i \neq j$.  This implies that $P(Y \vert Z=0, X=x_i)=P(Y \vert Z=0, X=x_j)$ and $P(Y  \vert Z=1, X=x_i)=P(Y \vert Z=1, X=x_j)$. We can write $p_i^k =p(Y_k=1 \vert X=x_i)=p(Y=1 \vert Z=k, X=x_i)$ for $k=0,1$ for simplicity. If we select the balancing score as a paired-value score, $b(x_i)=(p_i^0,p_i^1)$ then state space of $b(X)$ has a cardinality that is lesser than that of $X$ by  $1$ if  $b(x_i)=b(x_j)$,  for some $i \neq j$. We call $b(x)$ the potential outcome score, or simply the outcome score though it is a pair of scales unlike the propensity score. To see this, as previously if we define the random variable $X'$ from $X$ in the case of $i=1$ and $j=2$ (or in the reverse order) then we have  
\begin{align*}
 \sum_{i=1}^2 & p(Y=y \vert X=x_i,Z=1)p(X=x_i)  - \sum_{i=1}^2 p(Y=y \vert X=x_i,Z=0)p(X=x_i)  \\
   &= p(Y=y \vert (X=x_1 \textrm{ or } X=x_2),Z=1)  \sum_{i=1}^2 P(X=x_i) - p(Y=y \vert X=x_1 \textrm{ or } X=x_2),Z=0) \sum_{i=1}^2 P(X=x_i) \\
  &=  p(Y=y \vert (X=x_1 \textrm{ or } X=x_2),Z=1) p(X=x_1 \textrm{ or } X=x_2)  - P(Y=y \vert (X=x_1 \textrm{ or } X=x_2),Z=0) p(X=x_1 \textrm{ or } X=x_2)  \\
  &= p(Y=y \vert X'=x'_1,Z=1)p(X'=x'_1)  - p(Y=y \vert X'=x'_1,Z=0)p(X'=x'_1)
\end{align*}
Therefore we have from Eq. (\ref{tau})
\begin{align*}
\tau &= \sum_y y \sum_{i=1}^{n-1} p(Y=y \vert X'=x'_i,Z=1)p(X'=x'_i)  - \sum_y y \sum_{i=1}^{n-1} p(Y=y \vert X'=x'_i,Z=0)p(X'=x'_i) 
\end{align*}
That is, we get the same causal effect estimate using $X'$ as the conditioning variable. So, dimension reduction with the outcome score is also possible  as in the case of the propensity score.  And computational advantage of the outcome score is similar to that of the propensity score.  

Interpreting of the potential outcome score $b(x)=(p(Y_0=1 \vert X=x), p(Y_1=1 \vert X=x))$ should be  done in the line of doing that for the propensity score, i.e., in terms of mimicking randomization in observational studies. In a randomized study the probability distributions of the potential outcomes are assumed to be the same within the whole treatment group and control group separately irrespective of values of $X$ of subjects in them. That is, the outcome is not confounded. Therefore comparison of two groups as wholes is possible. But in an observational study, the potential outcomes are assumed to be dependent on the values of $X$ for each subject, therefore comparisons should be done subgroup-wise with respective to values of $X$ because within each strata of $X$ the potential outcome distributions are the same in treatment group and control group separately. So, if it is found that the potential outcomes distributions are the same for two values of $X$, say, $x_i$ and $x_j$, i.e., $b(x_i)=b(x_j)$, then it is not wrong to compare treatment and control subgroups corresponding to $X=x_i$ and $X=x_j$ after pooling these two difference configurations of $X$ together. Recall that when $P(Y_k \vert X=x_i)=P(Y_k \vert X=x_j)$ implies that they are equal to $P(Y_k \vert X=x_i \textrm{ or } X=x_j)$ for $k=0,1$. 

When only strict equalities are used the propensity score defines  conditioning variable that has state space with minimal cardinality  when using only $Z$ and $X$. And similarly outcome score gives that using  
$Y$, $Z$ and $X$. Since both scores behave independent of each other, it is possible to do dimension reduction in the conditioning variable $X$ in terms both of them. We state it in a theorem as follows since its practical value.

\textbf{Theorem 3}~\emph{When strictly equalities separately among the values of the propensity scores and the outcome scores are used the two scores  together can define a conditioning variable with a state space with minimal cardinality.}

Though we think in the line of comparing subgroups with respective to values of $X$, i.e., stratification based on confounding variable, unlike in the propensity score based stratification it is possible to stratify the treatment group and the control group separately in the case of using the potential outcome scores. This is evident from the causal effect estimation formula.  Subgroup-wise comparison of the treatment and the control groups  means algebraically taking the difference between the weighted averaged expected outcome in the treatment group and that of the control group where weights are marginal probabilities of $X$. And further two terms that are being weighted corresponds to two conditional probabilities $P(Y \vert X, Z=1)$ and $P(Y \vert X, Z=0)$ that are in fact outcome scores. Therefore we can have two separate new stratifications, one for the treatment group and the other for the control group. To see this simply, let us assume that we have four potential outcome scores such that they are component-wise equal as follows; without loss of generality, take  $b(x_1)$, $ b(x_2)$,  $b(x_3)$ and $b(x_4)$ such that $p_1^1=p_2^1$ and $p_3^0=p_4^0$ where $n>4$. Then, 
\begin{align*}
 \sum_{i=1}^n & p(Y=y \vert X=x_i,Z=1)p(X=x_i)  - \sum_{i=1}^n p(Y=y \vert X=x_i,Z=0)p(X=x_i)  \\
   &=  p(Y=y \vert (X=x_1 \textrm{ or } X=x_2),Z=1) p(X=x_1 \textrm{ or } X=x_2) + \sum_{i=3}^n p(Y=y \vert X=x_i,Z=1)p(X=x_i) \\
   & \qquad{}   - \sum_{i=1}^2 p(Y=y \vert X=x_i,Z=0) p(X=x_i) \\  
   &\qquad{}  - p(Y=y \vert ( X=x_3 \textrm{ or } X=x_4),Z=0) p(X=x_3 \textrm{ or } X=x_4)  - \sum_{i=5}^n p(Y=y \vert X=x_i,Z=0)p(X=x_i)   \\
  &=  \sum_{i=1}^{n-1} p(Y=y \vert X'=x'_i,Z=1)p(X'=x'_i)  - \sum_{i=1}^{n-1} p(Y=y \vert X''=x''_i,Z=0)p(X''=x''_i) 
\end{align*}
where $X'$ is defined as above using $p_1^1=p_2^1$ and similarly $X''$ is defined using $p_3^0=p_4^0$, both from $X$ with appropriate indexing. Therefore
\begin{align*}
\tau = & \sum_y y  \sum_{i=1}^{n-1} p(Y=y \vert X'=x'_i,Z=1)p(X'=x'_i) - \sum_y y \sum_{i=1}^{n-1} p(Y=y \vert X''=x''_i,Z=0)p(X''=x''_i)
\end{align*}
That is, two separate stratifications also work. This is true for any set of values of $X$ satisfying equalities among conditional probabilities  $p(Y=y \vert Z=1,X=x)$ and  $p(Y=y \vert Z=0,X=x)$ separately. However we are not getting a single  confounding variable with reduced state space from the initial confounding variable $X$ if used outcome score pairs refer to distinct configurations of $X$, for example, when $p^1_i=p^1_j$ for $i \neq j$ and  $p^0_k=p^0_l$ for $k \neq l$ but unordered pairs $(i,j)$ and $(k,l)$ differ.
\vspace{0.3cm}

\textbf{4. Balancing Scores and Causal Graphical Models}

One question that may be of interest is that if the propensity scores $L(X)=p(Z=1 \vert X)$ can be used in the graphical models. In literature (Pearl 2009) it is stated that its independence property, i.e., $X \perp Z \vert L(X)$ cannot be represented in graphs considering the case that $L(X)$ should be a child node (effect variable) of $X$ in the causal graph on $Z,Y$ and $X$ (see Fig. \ref{simple.bn}). And the conclusion was that $L(X)$ would not in general d-separate $Z$ from $X$. However from our above analysis of the propensity scores it is clear that when some equalities exist among values of $L(X)$ we are forming a new conditioning variable (confounder) $X'$ from $X$, that has a lesser cardinality in its state space  than that of $X$ and has well defined probability distribution. Note that we are in the case of joint probability distribution of $Z,Y$ and $X$ that exists. That is we transform $(Z,Y,X)$ contingency table to that of $(Z,Y,X')$. So, we can have a new graphical model on   $Z,Y$ and $X'$   from the old one that is on  $Z,Y$ and $X$ by replacing $X$ by $X'$. Both are equivalent in terms of yielding same causal effect estimates.   This new graph implicitly represents the variable $L(X)$ so there is no real need to have $L(X)$ explicitly in the graphical model on $Z,Y$ and $X$. Furthermore note that one can have causal relations $X \rightarrow X'$, $X' \rightarrow Z$, $X' \rightarrow Y$ and $Z \rightarrow Y$ in the graph. Then $X'$ d-separates $Z$ and $X$ so does $L(X)$ implicitly. But for the causal effect calculation with the propensity score above first causal arrow is not needed so the new graph can have only the second three.  So, one can simply replace $X$ by $X'$ in the causal graph on $Z$, $Y$ and $X$ to obtain the new one on $Z$, $Y$ and $X'$. In fact, one can have $X$ replaced by $L(X)$ to get the new causal model.

However propensity scores are not the only balancing score therefore one may be interested in having new causal graph with all possible balancing scores (at least with propensity scores and outcome scores). But when outcome scores are used as we have seen, one can sometimes have two distinct new confounder variables ($X'$ and $X''$) created  simultaneously. But together they do not represent a single variable. And also if one applies both the propensity scores and the outcome scores it may be difficult to represent all those operations with one single variable in addition to the fact that there is no real need that we have one single representation of all those operations that was done consecutively.   Therefore, as far as the causal effect estimation using graphical causal models are concerned inclusion of balancing scores in them gives no advantage.  
\vspace{0.3cm}

\textbf{4. Discussion}

We have presented a probabilistic analysis of balancing scores such as propensity scores and prognostic scores used in the potential outcome causal inference  framework. We could show how to derive more general prognostic scores that we call outcome scores.  Our analysis allows us to derive both the scores simultaneously and therefore it gives an understanding of using them together. One can use both scores effectively to obtain a new confounder conditioning set with minimal cardinality in its state state space from the original confounder conditioning set. 

In fact, the propensity scores are the most popular balancing scores among the applied statistics community and the other users of statistical theory and the prognostic scores are used relatively very rarely. The main reason for this imbalance may be due to the way in which the prognostic scores are derived in previous literature. But from our derivation it is clear that the prognostic (or outcome) scores have the same, if not very similar, potentials as propensity scores do for using in causal inference tasks. We hope this derivation will help to popularize prognostic scores among the users of statistics more and more. Finally as a future research, since propensity scores have many uses, we wish to look into other uses of outcome (prognostic) scores. 

\vspace{0.3cm}

\textbf{Acknowledgements}

This research is financed by the Swedish Research Council through the Swedish Initiative for Microdata Research in the Medical and Social Sciences (SIMSAM)  and the Swedish Research Council for Health, Working Life and Welfare (FORTE).

\textbf{References}

\hangafter=1
\setlength{\hangindent}{2em}
Cai, Z. and Kuroki, M. (2008). On Identifying Total Effects in the Presence of Latent Variables and Selection Bias. \emph{Proceedings of the Twenty-Fourth Conference Conference on Uncertainty in Artificial Intelligence, Helsinki, Finland (UAI-08)}, 62-69, Corvallis, Oregon, USA: AUAI Press. http://arxiv.org/abs/1206.3239

\hangafter=1
\setlength{\hangindent}{2em}
Glymour, C. (2013). Counterfactuals, graphical causal models, potential outocmes: Response to Lindquist and Sobel.  \emph{NeuroImage, 76}, 450-451.
http://dx.doi.org/10.1016/j.neuroimage.2011.07.071

\hangafter=1
\setlength{\hangindent}{2em}
Hansen, B. B. (2006). Bias Reduction in Observational Studies via Prognosis Scores. 
\newblock \emph{Technical Report, 441}, Statistics Department, University of Michigan, USA.

\hangafter=1
\setlength{\hangindent}{2em}
Hansen, B. B.  (2008). The Prognostic Analogue of the Propensity Score. \emph{Biometrika, 95}(2), 481-488. \linebreak 
http://dx.doi.org/10.1093/biomet/asn004

\hangafter=1
\setlength{\hangindent}{2em}
Holland, P. W. (1986). Statistics and Causal Inference. \emph{Journal of the American Statistical Association, 81}(396), 945-960. 

\hangafter=1
\setlength{\hangindent}{2em}
Lauritzen, S. L., and Richardson, T. S. (2002). 
\newblock Chain graph models and their causal interpretations.
\newblock  \emph{Journal of Royal Statistical Society, Series B, 64}(3), 321-361.

\hangafter=1
\setlength{\hangindent}{2em}
Martens E. P, Pestman, W. R., de Boer, A. , Belitser, S. V.  and  Klungel, O. H. (2008). Systematic differences in treatment effect estimates between propensity score methods and logistic regression. \emph{International Journal of Epidemiology, 37}(5), 1142-1147.
http://dx.doi.org/10.1093/ije/dyn079

\hangafter=1
\setlength{\hangindent}{2em}
Pearl, J.  (2009). \emph{Causality: Models, Reasoning, and Inference}, New York, USA: Cambridge University Press.

\hangafter=1
\setlength{\hangindent}{2em}
Pearl, J. (2009a). Letter to the Editor: Remarks on the method of propensity score. \emph{Statistics in Medicine, 28}(9), 1415-1416. http://dx.doi.org/10.1002/sim.3521

\hangafter=1
\setlength{\hangindent}{2em}
Rosenbaum P.R., Rubin D.B. (1983). The central role of the propensity score in observational studies for causal effect. \emph{Biometrika, 70}(1), 41-55.

\hangafter=1
\setlength{\hangindent}{2em}
Rubin, Donald(1974). Estimating Causal Effects of Treatments in Randomized and Nonrandomized Studies. \emph{Journal of Educational Psychology, 66}(5), 688-701.

\hangafter=1
\setlength{\hangindent}{2em}
Rubin, D. B. (2007). The design versus the analysis of observational studies for causal effects: parallels with the design of randomized trials. \emph{Statistics in Medicine, 26}(1), 20-36.  http://dx.doi.org/10.1002/sim.2739

\hangafter=1
\setlength{\hangindent}{2em}
Rubin, D. B. (2008). For objective causal inference, design trumps analysis. \emph{The Annals of Applied Statistics, 2}(3), 808-840.

\hangafter=1
\setlength{\hangindent}{2em}
Rubin, D. (2009). Author's Reply: Should observational studies be designed to allow lack of balance in covariate distributions across treatment groups?  \emph{Statistics in Medicine, 28}(9), 1420-1423. http://dx.doi.org/10.1002/sim.3565

\hangafter=1
\setlength{\hangindent}{2em}
Schneeweiss, S., Rassen, J. A., Glynn, R. J., Avorn, J., Mogun, H., and Brookhart, M. A. (2009). High-dimensional propensity score adjustment in studies of treatment effects using health care claims data.  \emph{Epidemiology, 20}(4), 512-522. http://dx.doi.org/10.1097/EDE.0b013e3181a663cc

\hangafter=1
\setlength{\hangindent}{2em}
Shah, B. R, Laupacis, A, Hux, J. E., Austin, P. C. (2005). Propensity score methods gave similar results to traditional regression modeling in observational studies: a systematic review. 
\newblock \emph{Journal of Clinical Epidemiology, 58}(6), 550-559. 
http://dx.doi.org/10.1016/j.jclinepi.2004.10.016

\hangafter=1
\setlength{\hangindent}{2em}
Shrier, I. (2008). Letter to the Editor. \emph{Statistics in Medicine, 27}(14), 2740-2741. http://dx.doi.org/10.1002/sim.3172

\hangafter=1
\setlength{\hangindent}{2em}
Sj\"{o}lander, A. (2009). Letter to the Editor. \emph{Statistics in Medicine, 28}(9), 1416-1420.
http://dx.doi.org/10.1002/sim.3532

\hangafter=1
\setlength{\hangindent}{2em}
Spirtes, P., Glymour, C. and Scheines, R. \emph{Causation, Prediction, and Search}, 2nd edition, (2001). USA: MIT Press.

\hangafter=1
\setlength{\hangindent}{2em}
Stuart, A. E., Brian, K. L. and Leacy, F. P. (2013). Prognostic score-based balance measures can be a useful diagnostic for propensity score methods in compartaive effectiveness research. \emph{Journal of Clinical Epidemiology, 66}(8 Suppl), S84-S90. http://dx.doi.org/10.1016/j.jclinepi.2013.01.013

\hangafter=1
\setlength{\hangindent}{2em}
St\"urmer T, Joshi, M, Glynn R.J, Avorn J, Rothman KJ, Schneeweiss S. (2006). A review of the application of propensity score methods yielded increasing use, advantages in specific settings, but not substantially different estimates compared with conventional multivariable methods. \emph{Journal of  Clinical Epidemiology, 59}(5), 437-447.

\hangafter=1
\setlength{\hangindent}{2em}
Vansteelandt, S., Bekaert, M., and Claeskens, G. (2010).  On Model Selection and Model Misspecification in Causal Inference. \emph{Statistical Methods in Medical Research, 21}(1), 7-30. http://dx.doi.org/10.1177/0962280210387717

\hangafter=1
\setlength{\hangindent}{2em}
Vansteelandt, S. (2012). Discussions. \emph{Biometrics, 68}(3), 675-678. http://dx.doi.org/10.1111/j.1541-0420.2011.01734.x

\hangafter=1
\setlength{\hangindent}{2em}
Wijayatunga, P. (2014). Causal Effect Estimation Methods. \emph{Journal of Statistical and Econometric Methods, 3}(2), 153-170. http://www.scienpress.com/Upload/JSEM/Vol{\%}203{\_}2{\_}9.pdf

\hangafter=1
\setlength{\hangindent}{2em} 
Williamson, E. J., Morley, R., Lucas, A. and Carpenter, J. R. (2012).  Variance estimation for stratified propensity score estimators. \emph{Statistics in Medicine, 31}(15), 1617-1632. 
http://dx.doi.org/10.1002/sim.4504

\hangafter=1
\setlength{\hangindent}{2em}
Zigler, C. M.,1, Watts,  K., Yeh, R. W., Wang, Y. , Coull, B. A., and Dominici, F. (2013). Model Feedback in Bayesian Propensity Score Estimation. \emph{Biometrics, 69}(1), 263-273. http://dx.doi.org/10.1111/j.1541-0420.2012.01830.x

\vspace{0.5cm}
\textbf{Copyrights}

Copyright for this article is retained by the author(s), with first publication rights granted to the journal.

This is an open-access article distributed under the terms and conditions of the Creative Commons Attribution license (http://creativecommons.org/licenses/by/3.0/).

\end{document}